# Variance Analysis of Multiple Importance Sampling Schemes

Rahul Mukerjee　　　　　　　　and　　　　Víctor Elvira
Indian Institute of Management Calcutta　　　　School of Mathematics
Joka, Diamond Harbour Road　　　　　　　　　University of Edinburgh
Kolkata 700104, India　　　　　　　　　　　　Edinburgh EH9 3FD, United Kingdom
E-mail: rmuk0902@gmail.com　　　　　　　　　E-mail: victor.elvira@ed.ac.uk

*Abstract*: Multiple importance sampling (MIS) is an increasingly used methodology where several proposal densities are used to approximate integrals, generally involving target probability density functions. The use of several proposals allows for a large variety of sampling and weighting schemes. Then, the practitioner must choose a given scheme, i.e., sampling mechanism and weighting function. A variance analysis has been proposed in Elvira et al (2019, *Statistical Science* **34**, 129-155), showing the superiority of the balanced heuristic estimator with respect to other competing schemes in some scenarios. However, some of their results are valid only for two proposals. In this paper, we extend and generalize these results, providing novel proofs that allow to determine the variance relations among MIS schemes.

## 1. Introduction

Importance sampling (IS) methods are emerging Monte Carlo techniques for the approximation of integrals, generally involving probability density functions (Robert and Casella, 2004; Liu, 2004; Tokdar and Kass, 2010; Chatterjee and Diaconis, 2018; Elvira and Martino, 2021). In the standard IS algorithm, the samples are simulated from a proposal density, and each sample receives an importance weight. This weight allows correcting the mismatch between the proposal density and the target density in the integrand. The efficiency of IS methods, generally measured in terms of the variance of the estimators, depends on the mismatch between the proposal and target densities. For instance, it is well known that the variance of the standard estimator, also called unnormalized IS (UIS) estimator, is strongly connected to the $\chi^2$-divergence between target and proposal (Ryu and Boyd, 2014; Agapiou et al., 2017).

The multiple IS (MIS) paradigm (Owen and Zhou, 2000; Sbert et al., 2018; Elvira et al., 2019), i.e., the use of several proposal densities, has gained a considerable attention in recent years, primarily due to two reasons. First, mixture distributions present an interesting flexibility for the approximation of a wide range of distributions. Second, it is known to be a difficult problem to find one single proposal that minimizes its mismatch with respect to the target. For that reason, many adaptive IS algorithms have been proposed in the last few years, e.g., Martino et al. (2014, 2017), Fasiolo et al. (2018), Elvira and Chouzenoux (2022); see Bugallo et al. (2017) for a review.

The use of several proposal densities offers methodological opportunities since the possible sampling and weighting schemes are almost unlimited; Elvira et al. (2015). Recently, Elvira et al. (2019), proposed a general framework of MIS methods, and compared existing and new methods in terms of variance of the UIS estimator. With reference to a set of $N$ proposal densities, they considered three selection strategies, in conjunction with five choices of weighting functions, and arrived at six distinct MIS schemes, namely, N1, N2, N3, for selection without replacement, and R1, R2, R3, for selection with replacement; a description of these schemes requires more notation and appears later in this paper.

In Section 6 of their paper, Elvira et al. (2019), focused on a variance analysis of the aforesaid six schemes. One of their main results, Theorem 6.2, pertained to $N = 2$ proposal densities, and the case of general $N$ was posed by them as an open problem. The present work addresses this problem and strengthens the theorem in several directions. First, we show that the scheme N2 is superior to N1 (no larger variance) for any target density, test function, and set of $N$ ($\geq 2$) proposals. Then, we also show that N2 is inferior to N3 (no smaller variance) for any target density, test function, and set of $N$ ($\geq 2$) proposals. Finally, we show that R2 is superior to R1 or N1, again with the same level of generality. The cornerstone of our results on N2 is an expression for the corresponding variance which is derived here for general $N$ via a sequence of arguments, some combinatorial.

**2. Selection strategies**

For simplicity, we write $u(x)$ for what Elvira et al. (2019) call $\pi(x)g(x)/Z$ in their notation, with $\pi(x)$ proportional to the target density, $g(x)$ denoting is any square integrable function with respect to $\pi(x)$, and $Z = \int \pi(x)dx$, the range of integration, say $\mathcal{X}$, being suppressed to simplify notation. Interest lies in estimating $I = \int u(x)dx$. Let $q_1(x),....,q_N(x)$ be $N$ proposal probability density functions over $\mathcal{X}$. In MIS, the integral $I$ is estimated on the basis of samples $X_1,...., X_N$ where $X_n$ is randomly generated from $q_{j_n}(x)$ ($n = 1,…, N$). As in Elvira et al. (2019), we consider three possible ways for selecting the indices $j_1,..., j_N$:

(a) Random selection with replacement (RSWR): $j_1,..., j_N$ are independent, each uniformly distributed over $\{1,…,N\}$.

(b) Random selection without replacement (RSWoR): $(j_1,..., j_N)$ is a random permutation of $(1,…, N)$, all such $N!$ permutations being equally likely,

(c) Deterministic selection without replacement (DSWoR): $(j_1,..., j_N) = (1,…, N)$.

Throughout, given $(j_1,..., j_N)$, the random variables $X_1,...., X_N$ are supposed to be independent. Moreover, all expectations and variances in what follows are assumed to exist.

## 3. Index selection without replacement

This section dwells on RSWoR and DSWoR. In this setup, the schemes N1, N2 and N3 in Elvira et al. (2019) correspond to three estimators of $I$, namely,

$$\hat{I}_{N1} = (1/N)\Sigma_{n=1}^{N} u(X_n)/q_{j_n}(X_n), \tag{1}$$

$$\hat{I}_{N2} = (1/N)\Sigma_{n=1}^{N} u(X_n)/\lambda_{nj}(X_n), \tag{2}$$

$$\hat{I}_{N3} = (1/N)\Sigma_{n=1}^{N} u(X_n)/\psi(X_n), \tag{3}$$

where

$$\lambda_{nj}(x) = \{1/(N-n+1)\}\Sigma_{k=n}^{N} q_{j_k}(x), \quad n = 1,\ldots, N. \tag{4}$$

with $j = (j_1,\ldots, j_N)$, and

$$\psi(x) = (1/N)\Sigma_{n=1}^{N} q_n(x). \tag{5}$$

Note that $\lambda_{1j}(x) = \psi(x)$, as $(j_1,\ldots, j_N)$ is a permutation of $(1,\ldots, N)$. Elvira et al. (2019) considered $\hat{I}_{N1}$ and $\hat{I}_{N3}$ for both RSWoR and DSWoR, and $\hat{I}_{N2}$ only for RSWoR. They obtained the following result on $\hat{I}_{N1}$ and $\hat{I}_{N3}$.

**Proposition 1**. *Under both* RSWoR *and* DSWoR,

$$E(\hat{I}_{N1}) = E(\hat{I}_{N3}) = I,$$

$$\text{var}(\hat{I}_{N1}) = (1/N^2)\Sigma_{n=1}^{N}\int\{u^2(x)/q_n(x)\}dx - (1/N)I^2,$$

$$\text{var}(\hat{I}_{N3}) = (1/N)\int\{u^2(x)/\psi(x)\}dx - (1/N^2)\Sigma_{n=1}^{N}[\int\{u(x)/\psi(x)\}q_n(x)dx]^2.$$

A key point in the proof of the variance formulae in Proposition 1 under RSWoR is that the none of the conditional expectations $E(\hat{I}_{N1} \mid j)$ and $E(\hat{I}_{N3} \mid j)$, where $j = (j_1,\ldots, j_N)$ is a permutation of $(1,\ldots, N)$, depends on $j$. Hence, $\text{var}(\hat{I}_{N1})$ and $\text{var}(\hat{I}_{N3})$ can simply be obtained as the expectations of the corresponding conditional variances $\text{var}(\hat{I}_{N1} \mid j)$ and $\text{var}(\hat{I}_{N3} \mid j)$. As seen below, this simplifying feature is not shared by $\hat{I}_{N2}$, i.e., under RSWoR, the conditional expectation $E(\hat{I}_{N2} \mid j)$ can depend on $j$, and hence $\text{var}\{E(\hat{I}_{N2} \mid j)\}$ can be positive. Therefore, one needs to obtain $\text{var}(\hat{I}_{N2})$ as

$$\text{var}(\hat{I}_{N2}) = E\{\text{var}(\hat{I}_{N2} \mid j)\} + \text{var}\{E(\hat{I}_{N2} \mid j)\}. \tag{6}$$

This point was missed in Elvira et al. (2019) who obtained $\text{var}(\hat{I}_{N2})$ only as $E\{\text{var}(\hat{I}_{N2} \mid j)\}$ [cf. their equation (D.4)], ignoring the second term in (6). The resulting expression for $\text{var}(\hat{I}_{N2})$ is, therefore, incomplete. However, as we show here, the inequality $\text{var}(\hat{I}_{N1}) \geq \text{var}(\hat{I}_{N2}) \geq \text{var}(\hat{I}_{N3})$ that they ob-tained for $N = 2$ in their Theorem 6.2 still remains valid. Even more, we prove the truth of this inequality for general $N$ and Theorem 1 below, showing the correct form of $\text{var}(\hat{I}_{N2})$, is crucial for this purpose.

For completeness, before presenting Theorem 1, we indicate another fact about $\hat{I}_{N2}$. Elvira et al. (2019) considered this estimator only for RSWoR and not for DSWoR, even though it is well defined in the latter case with $(j_1,\ldots, j_N) = (1,\ldots, N)$. Doing so is justified in view of the next example which shows that $\hat{I}_{N2}$ is in general biased for $I$ under DSWoR.

**Example 1**. Let $\mathcal{X} = [0, 2]$, $N = 2$, and

$$q_1(x) = \begin{cases} b & \text{if } 0 \leq x \leq 1 \\ 1-b & \text{if } 1 < x \leq 2 \end{cases}, \quad q_2(x) = \begin{cases} 1-b & \text{if } 0 \leq x \leq 1 \\ b & \text{if } 1 < x \leq 2 \end{cases}.$$

where $b$ is a constant satisfying $0 < b < 1$, $b \neq 1/2$. Note that $\{q_1(x) + q_2(x)\}/2 = 1/2$ over $\mathcal{X}$ and that $I = I_1 + I_2$, where $I_1 = \int_0^1 u(x)dx$ and $I_2 = \int_1^2 u(x)dx$. Hence by (2) and (4), with $(j_1, j_2) = (1, 2)$,

$$E(\hat{I}_{N2}) = (1/2)[\int_0^2 \{u(x)/(1/2)\}q_1(x)dx + \int_0^2 \{u(x)/q_2(x)\}q_2(x)dx]$$

$$= \int_0^2 u(x)q_1(x)dx + (1/2)I = bI_1 + (1-b)I_2 + (1/2)I,$$

Since $b \neq 1/2$, the above differs from $I$ whenever $I_1 \neq I_2$. □

We now present Theorem 1. Its proof, involving a sequence of steps, appears in the appendix. Let us establish some useful notation. For $n = 1,\ldots, N$, let $A_n$ be the class of subsets of $\{1,\ldots, N\}$ having cardinality $N - n + 1$. Clearly, there are $\binom{N}{n-1}$ sets in $A_n$. For any set $a$ in $A_n$, let

$$\xi_a(x) = \{1/(N-n+1)\}\Sigma_{\gamma \in a} q_\gamma(x). \tag{7}$$

**Theorem 1**. *Under* RSWoR, $E(\hat{I}_{N2}) = I$, *and*

$$\text{var}(\hat{I}_{N2}) = (1/N^2)\Sigma_{n=1}^N \{1/\binom{N}{n-1}\}\Sigma_{a \in A_n} \int \{u^2(x)/\xi_a(x)\}dx - (1/N)I^2.$$

In particular, for $N = 2$, using (5) and (7), Theorem 1 yields

$$\text{var}(\hat{I}_{N2}) = (1/4)\int \{u^2(x)/\psi(x)\}dx + (1/8)\Sigma_{n=1}^2 \int \{u^2(x)/q_n(x)\}dx - (1/2)I^2, \tag{8}$$

while proceeding as in the appendix,

$$E\{\text{var}(\hat{I}_{N2} \mid j)\} = (1/4)\int \{u^2(x)/\psi(x)\}dx + (1/8)\Sigma_{n=1}^2 \int \{u^2(x)/q_n(x)\}dx$$

$$- (1/8)\Sigma_{n=1}^2 [\int \{u(x)/\psi(x)\}q_n(x)dx]^2 - (1/4)I^2, \tag{9}$$

$$\text{var}\{E(\hat{I}_{N2} \mid j)\} = (1/8)\Sigma_{n=1}^2 [\int \{u(x)/\psi(x)\}q_n(x)dx]^2 - (1/4)I^2. \tag{10}$$

Equations (8)-(10) agree with (6). On the other hand, Elvira et al. (2019) reported $\text{var}(\hat{I}_{N2})$ as only (9), ignoring (10) which, by (5), is positive unless $\int \{u(x)/\psi(x)\}q_1(x)\mathrm{d}x = \int \{u(x)/\psi(x)\}q_2(x)\mathrm{d}x$.

As a comparison of (8) and (9) suggests, $\text{var}(\hat{I}_{N2})$ is appreciably simpler than $E\{\text{var}(\hat{I}_{N2} \mid j)\}$. This has implications regarding the variance inequality $\text{var}(\hat{I}_{N1}) \geq \text{var}(\hat{I}_{N2}) \geq \text{var}(\hat{I}_{N3})$, that Elvira et al. (2019) established only for $N = 2$, with $\text{var}(\hat{I}_{N2})$ taken as the expression in (9). Reassuringly, we find that this inequality remains valid even if one works with the correct variance formula in Theorem 1. Indeed, the relative simplicity of this variance formula allows us to prove the inequality for general $N$, thus strengthening the findings in Elvira et al. (2019). This is done in Theorem 2 below which is proved in the appendix.

**Theorem 2**. *For any set of proposal densities $q_1(x),...,q_N(x)$,*

$$\text{var}(\hat{I}_{N1}) \geq \text{var}(\hat{I}_{N2}) \geq \text{var}(\hat{I}_{N3}).$$

**4. Random index selection with replacement (RSWR)**

We now turn to RSWR. In this setup, the schemes R1, R2 and R3 in Elvira et al. (2019) correspond to three estimators of $I$, namely,

$$\hat{I}_{R1} = (1/N)\Sigma_{n=1}^{N} u(X_n)/q_{j_n}(X_n), \qquad \hat{I}_{R2} = (1/N)\Sigma_{n=1}^{N} u(X_n)/\bar{q}_j(X_n),$$

$$\hat{I}_{R3} = (1/N)\Sigma_{n=1}^{N} u(X_n)/\psi(X_n),$$

where, for any $j = (j_1,...,j_N)$,

$$\bar{q}_j(x) = (1/N)\Sigma_{n=1}^{N} q_{j_n}(x). \tag{11}$$

Recalling (1) and (3), $\hat{I}_{R1}$ and $\hat{I}_{R3}$ are formally identical to $\hat{I}_{N1}$ and $\hat{I}_{N3}$, respectively, but they refer to index selection with, rather than without, replacement. This distinction is reflected in the following result, obtained by Elvira et al. (2019), which shows that $\text{var}(\hat{I}_{R3})$ is different from $\text{var}(\hat{I}_{N3})$ in Proposition 1, even though $\text{var}(\hat{I}_{R1})$ equals $\text{var}(\hat{I}_{N1})$.

**Proposition 2**. *Under RSWR,*

$$E(\hat{I}_{R1}) = E(\hat{I}_{R2}) = E(\hat{I}_{R3}) = I,$$

$$\text{var}(\hat{I}_{R1}) = (1/N^2)\Sigma_{n=1}^{N} \int \{u^2(x)/q_n(x)\}\mathrm{d}x - (1/N)I^2,$$

$$\text{var}(\hat{I}_{R2}) = (1/N^N)\Sigma_j W(j),$$

$$\text{var}(\hat{I}_{R3}) = (1/N)\int \{u^2(x)/\psi(x)\}\mathrm{d}x - (1/N)I^2,$$

*where $\Sigma_j$ denotes sum over $j_1,...,j_N$, each with range of summation $\{1,...,N\}$, and*

$$W(j) = (1/N)\int\{u^2(x)/\overline{q}_j(x)\}dx - (1/N^2)\Sigma_{n=1}^{N}[\int\{u(x)/\overline{q}_j(x)\}q_{j_n}(x)dx]^2.$$

For general $N$, Elvira et al. (2019) proved that $\text{var}(\hat{I}_{R1}) = \text{var}(\hat{I}_{N1}) \geq \text{var}(\hat{I}_{R3}) \geq \text{var}(\hat{I}_{N3})$. From Theorems 1 and 2, Proposition 2 and (A.11) in the appendix, we also have $\text{var}(\hat{I}_{N1}) \geq \text{var}(\hat{I}_{N2}) \geq \text{var}(\hat{I}_{R3})$. Combining these, for general $N$,

$$\text{var}(\hat{I}_{R1}) = \text{var}(\hat{I}_{N1}) \geq \text{var}(\hat{I}_{N2}) \geq \text{var}(\hat{I}_{R3}) \geq \text{var}(\hat{I}_{N3}).$$

Turning to $\text{var}(\hat{I}_{R2})$, Elvira et al. (2019) showed that for $N = 2$,

$$\text{var}(\hat{I}_{R1}) \geq \text{var}(\hat{I}_{R2}) \geq \text{var}(\hat{I}_{N3}). \tag{12}$$

Our next result, proved in the appendix, establishes the first inequality in (12) for general $N$.

**Theorem 3**. *For any set of proposal densities $q_1(x),...,q_N(x)$,*

$$\text{var}(\hat{I}_{R1}) \geq \text{var}(\hat{I}_{R2}).$$

Computations, including those with piecewise constant proposal densities as in Example 1, suggest that the second inequality in (12) also holds for general $N$, though a proof of this remains intractable at this stage. From Propositions 1, 2 and (5), $W(j)$ equals $\text{var}(\hat{I}_{N3})$ if $j_1,...,j_N$ are distinct. Hence some generalized version of an inequality in Hoeffding (1963, Section 6) on expectations under simple random sampling with and without replacement may be useful in proving this second inequality for general $N$. We conclude with the hope that the present endeavor will generate further interest in this and related problems.

**Appendix: Proofs of theorems**

*Proof of Theorem 1*. The result on expectation is in agreement with Elvira et al. (2019), but we prove it here again in a unified framework that works for both expectation and variance, our main interest being in the latter.

Step 1 (conditional expectation and variance) For any $j = (j_1,..., j_N)$ which is a permutation of $(1,..., N)$, recalling the independence of $X_1,..., X_N$ given $(j_1,..., j_N)$, we get from (2) that

$$E(\hat{I}_{N2} \mid j) = (1/N)\Sigma_{n=1}^{N}\mu_n(j), \quad \text{var}(\hat{I}_{N2} \mid j) = (1/N^2)\Sigma_{n=1}^{N}V_n(j), \tag{A.1}$$

where, for $n = 1,..., N$,

$$\mu_n(j) = \int\{u(x)/\lambda_{nj}(x)\}q_{j_n}(x)dx, \quad V_n(j) = \int\{u(x)/\lambda_{nj}(x)\}^2 q_{j_n}(x)dx - \{\mu_n(j)\}^2. \tag{A.2}$$

Step 2 (some notation and useful facts) Let $n \geq 2$. For any permutation $j = (j_1,..., j_N)$ of $(1,..., N)$, let $j^* = (j_1,..., j_{n-1})$ and $\tilde{j} = (j_n,..., j_N)$. Write $B_n$ for the set of all $N!/(N - n + 1)!$ possibilities for $j^*$. Given any $j^* \in B_n$, let $C(j^*)$ be the set of all $(N - n + 1)!$ possibilities for $\tilde{j}$, and $a(j^*)$ denote the complement of $\{j_1,..., j_{n-1}\}$ in $\{1,..., N\}$. Thus, $a(j^*)$ is the unordered set corresponding to the ordered tuple $\tilde{j}$, so that $C(j^*)$ consists of all permutations of members of $a(j^*)$. Note that $a(j^*)$ belongs to the class of sets $A_n$ in Theorem 2. Hence, given $j^*$, by (4) and (7),

$$\lambda_{nj}(x) = \{1/(N-n+1)\}\Sigma_{\gamma \in a(j^*)} q_\gamma(x) = \xi_{a(j^*)}(x), \quad \text{(A.3)}$$

$$\Sigma_{\tilde{j} \in C(j^*)} q_{j_n}(x) = (N-n)!\Sigma_{\gamma \in a(j^*)} q_\gamma(x) = (N-n+1)!\xi_{a(j^*)}(x), \quad \text{(A.4)}$$

the first identity in (A.4) being a consequence of the fact that for every $\gamma \in a(j^*)$, there are $(N-n)!$ members of $C(j^*)$ for which $j_n$ equals $\gamma$.

<u>Step 3</u> (unconditional expectation) For $n \geq 2$ and each $j^* \in B_n$, by (A.2)-(A.4),

$$E\{\mu_n(j) | j^*\} = \{(1/(N-n+1)!\}\Sigma_{\tilde{j} \in C(j^*)} \int \{u(x)/\xi_{a(j^*)}(x)\} q_{j_n}(x) dx$$

$$= \int \{u(x)/\xi_{a(j^*)}(x)\} \xi_{a(j^*)}(x) dx = I. \quad \text{(A.5)}$$

The constancy of the above over $j^*$ implies that

$$E\{\mu_n(j)\} = I. \quad \text{(A.6)}$$

Because $\lambda_{1j}(x) = \psi(x)$, it is readily seen from (A.2) and (5) that (A.6) holds for $n = 1$ as well. Thus (A.6) holds for every $n = 1,..., N$, and by (A.1), $E(\hat{I}_{N2}) = I$.

<u>Step 4</u> (expectation of conditional variance) For $n \geq 2$ and each $j^* \in B_n$, analogously to (A.5),

$$E[\int \{u(x)/\lambda_{nj}(x)\}^2 q_{j_n}(x) dx | j^*] = \int \{u^2(x)/\xi_{a(j^*)}(x)\} dx.$$

Taking expectation of the above over $j^*$,

$$E[\int \{u(x)/\lambda_{nj}(x)\}^2 q_{j_n}(x) dx]$$

$$= \{(N-n+1)!/N!\}\Sigma_{j^* \in B_n} \int \{u^2(x)/\xi_{a(j^*)}(x)\} dx$$

$$= \{(N-n+1)!(n-1)!/N!\}\Sigma_{a \in A_n} \int \{u^2(x)/\xi_a(x)\} dx$$

$$= \{1/\binom{N}{n-1}\}\Sigma_{a \in A_n} \int \{u^2(x)/\xi_a(x)\} dx. \quad \text{(A.7)}$$

Here we invoke that fact that each $a \in A_n$, there are $(n-1)!$ choices of $j^*$ in $B_n$ such that $a(j^*) = a$.

For $n = 1$, $A_n$ has a single member $a = \{1,..., N\}$, and by (5) and (7), $\xi_a(x) = \psi(x)$ for this $a$. Since $\lambda_{1j}(x) = \psi(x)$, it is clear that (A.7) holds for $n = 1$ as well. Thus, by (A.1) and (A.2),

$$E\{\text{var}(\hat{I}_{N2} | j)\}$$

$$= (1/N^2)\Sigma_{n=1}^{N}\{1/\binom{N}{n-1}\}\Sigma_{a \in A_n} \int \{u^2(x)/\xi_a(x)\} dx - (1/N^2)\Sigma_{n=1}^{N} E[\{\mu_n(j)\}^2]. \quad \text{(A.8)}$$

As seen later, there is no need to find the second term in (A.8) explicitly.

<u>Step 5</u> (vanishing covariances) We now show that

$$\text{cov}\{\mu_m(j), \mu_n(j)\} = 0, \quad m, n = 1,\ldots, N, \quad m \neq n. \tag{A.9}$$

Without loss of generality, let $m < n$, so that $n \geq 2$. Given any $j^* \in B_n$, then by (4) and analogously to (A.3), $\lambda_{mj}(x)$ depends on $j$ only through $j^*$, so that by (A.2), the same happens for $\mu_m(j)$ as well. Hence, writing $\mu_m(j) = \beta_m(j^*)$, say,

$$E[\mu_m(j)\{\mu_n(j) - I\} \mid j^*] = \beta_m(j^*)E[\{\mu_n(j) \mid j^*\} - I] = 0,$$

by (A.5). Since the above holds for every $j^*$, the truth of (A.9) follows from (A.6).

Step 6 (unconditional variance) By (A.1), (A.6) and (A.9),

$$\text{var}\{E(\hat{I}_{N2} \mid j)\} = (1/N^2)\Sigma_{n=1}^{N}\text{var}\{\mu_n(j)\} = (1/N^2)\Sigma_{n=1}^{N}E[\{\mu_n(j)\}^2] - (1/N)I^2.$$

Adding the above to $E\{\text{var}(\hat{I}_{N2} \mid j)\}$ as given by (A.8), the result on $\text{var}(\hat{I}_{N2})$ follows. □

*Proof of Theorem 2.* For $n = 1,\ldots, N$, let $\bar{\xi}_n(x)$ be the arithmetic mean of $\xi_a(x)$, over $a \in A_n$. By (7),

$$\bar{\xi}_n(x) = \{1/(N-n+1)\}\{1/\binom{N}{n-1}\}\Sigma_{a \in A_n}\Sigma_{\gamma \in a} q_\gamma(x).$$

Since for each $s = 1,\ldots, N$, among the $\binom{N}{n-1}$ sets in $A_n$, there are $\binom{N-1}{N-n}$ that contain $s$, recalling (5),

$$\bar{\xi}_n(x) = \{1/(N-n+1)\}\{\binom{N-1}{N-n}/\binom{N}{n-1}\}\Sigma_{s=1}^{N} q_s(x) = \psi(x). \tag{A.10}$$

Therefore, invoking the arithmetic mean harmonic mean inequality,

$$(1/N^2)\Sigma_{n=1}^{N}\{1/\binom{N}{n-1}\}\Sigma_{a \in A_n}\int\{u^2(x)/\xi_a(x)\}dx$$

$$\geq (1/N^2)\Sigma_{n=1}^{N}\int\{u^2(x)/\bar{\xi}_n(x)\}dx = (1/N)\int\{u^2(x)/\psi(x)\}dx. \tag{A.11}$$

Also, by (5),

$$(1/N^2)\Sigma_{n=1}^{N}[\int\{u(x)/\psi(x)\}q_n(x)dx]^2$$

$$\geq (1/N)[\int\{u(x)/\psi(x)\}\{(1/N)\Sigma_{n=1}^{N}q_n(x)\}dx]^2 = (1/N)I^2. \tag{A.12}$$

From (A.11), (A.12), Proposition 1 and Theorem 1, $\text{var}(\hat{I}_{N2}) \geq \text{var}(\hat{I}_{N3})$.

Next, for $n = 1,\ldots, N$ and $a \in A_n$, invoking the arithmetic mean harmonic mean inequality again, by (7),

$$1/\xi_a(x) \leq \{1/(N-n+1)\}\Sigma_{\gamma \in a}\{1/q_\gamma(x)\}.$$

So, using the same arguments as in (A.10),

$$\{1/\binom{N}{n-1}\}\Sigma_{a \in A_n}\{1/\xi_a(x)\} \leq \{1/(N-n+1)\}\{1/\binom{N}{n-1}\}\Sigma_{a \in A_n}\Sigma_{\gamma \in a}\{1/q_\gamma(x)\}$$

$$= \{1/(N-n+1)\}\{\binom{N-1}{N-n}/\binom{N}{n-1}\}\Sigma_{s=1}^{N}\{1/q_s(x)\} = (1/N)\Sigma_{n=1}^{N}\{1/q_n(x)\}.$$

As a result, the first term of $\text{var}(\hat{I}_{N2})$ in Theorem 1 does not exceed the first term of $\text{var}(\hat{I}_{N1})$ in Proposition 1. Since these two variances differ only in their first term, $\text{var}(\hat{I}_{N2}) \leq \text{var}(\hat{I}_{N1})$. □

***Proof of Theorem 3***. The proof is similar to that of Theorem 2. For every $j = (j_1,...,j_N)$, by (11) and the arithmetic mean and harmonic mean inequality,

$$1/\bar{q}_j(x) \leq (1/N)\Sigma_{n=1}^{N}\{1/q_{j_n}(x)\}$$

while analogously to (A.12),

$$(1/N^2)\Sigma_{n=1}^{N}[\int\{u(x)/\bar{q}_j(x)\}q_{j_n}(x)\mathrm{d}x]^2 \geq (1/N)I^2.$$

so that $W(j)$ in Proposition 2 satisfies

$$W(j) \leq (1/N^2)\Sigma_{n=1}^{N}\int\{u^2(x)/q_{j_n}(x)\}\mathrm{d}x - (1/N)I^2$$

Because

$$(1/N^N)\Sigma_j\Sigma_{n=1}^{N}\{1/q_{j_n}(x)\} = (1/N^N)\Sigma_{n=1}^{N}\Sigma_j\{1/q_{j_n}(x)\} = \Sigma_{n=1}^{N}\{1/q_n(x)\},$$

it follows from Propositions 1 and 2 that

$$\mathrm{var}(\hat{I}_{R2}) = (1/N^N)\Sigma_j W(j) \leq (1/N^2)\Sigma_{n=1}^{N}\int\{u^2(x)/q_n(x)\}\mathrm{d}x - (1/N)I^2 = \mathrm{var}(\hat{I}_{R1}). \quad \square$$

**Acknowledgement**: The work of RM was supported by a grant from the Science and Engineering Research Board, Government of India. The work of VE was funded by a Leverhulme Research Fellowship.